\pdfoutput=1
\documentclass[12pt]{article}%
\usepackage{amsmath}
\usepackage{amsfonts}
\usepackage{amssymb}
\usepackage{graphicx}%
\providecommand{\U}[1]{\protect\rule{.1in}{.1in}}
\newtheorem{theorem}{Theorem}

\newtheorem{algorithm}[theorem]{Algorithm}

\newtheorem{case}[theorem]{Case}
\newtheorem{claim}[theorem]{Claim}

\newtheorem{conjecture}[theorem]{Conjecture}

\newtheorem{problem}[theorem]{Problem}

\begin{document}

\title{\textbf{Can Linear Superiorization Be Useful for Linear Optimization
Problems?}}
\author{Yair Censor\\Department of Mathematics\\University of Haifa\\Mt.\ Carmel, Haifa 3498838, Israel\\(yair@math.haifa.ac.il) }
\date{March 17, 2016. Revised: July 19, 2016. Revised: September 27, 2016.}
\maketitle

\begin{abstract}
Linear superiorization considers linear programming problems but instead of
attempting to solve them with linear optimization methods it employs
perturbation resilient feasibility-seeking algorithms and steers them toward
reduced (not necessarily minimal) target function values. The two questions
that we set out to explore experimentally are (i) Does linear superiorization
provide a feasible point whose linear target function value is lower than that
obtained by running the same feasibility-seeking algorithm without
superiorization under identical conditions? and (ii) How does linear
superiorization fare in comparison with the Simplex method for solving linear
programming problems? Based on our computational experiments presented here,
the answers to these two questions are: \textquotedblleft
yes\textquotedblright\ and \textquotedblleft very well\textquotedblright, respectively.

\end{abstract}


\section{Introduction\label{sec:intro}}

In this paper we propose the linear superiorization method as a tool for
handling linear programming problems. The linear superiorization method is not
guaranteed to find a minimum point of the linear optimization problem but it
steers the linear feasibility-seeking algorithm that it employs toward points
with reduced target function value. This task is not identical with that of
finding a minimizer to the linear programming problem but for huge sized
problems, it puts in the hands of the end-user a viable alternative to the
Simplex method of linear programming, against which we compared it here.

The paper relies on previous theoretical work about superiorization included
in the papers cited in the sequel, in particular \cite{cz14-feje}. Our working
tools here are only experimental computations. In spite of this, we are not
discussing computational issues per se but using computations as a tool in an
exploratory practical validation, style \textquotedblleft proof of
concept\textquotedblright\footnote{\textquotedblleft A proof of concept (POC)
or a proof of principle is a realization of a certain method or idea to
demonstrate its feasibility, or a demonstration in principle, whose purpose is
to verify that some concept or theory has the potential of being used. A proof
of concept is usually small and may or may not be complete\textquotedblright.
(https://en.wikipedia.org/wiki/Proof\_of\_concept).}.

\textbf{What is superiorization?} Many constrained optimization methods are
based on methods for unconstrained optimization that are adapted to deal with
constraints. Such is, for example, the class of projected gradient methods
wherein the unconstrained minimization inner step \textquotedblleft
leads\textquotedblright\ the process and a projection onto the whole
constraints set (the feasible set) is performed after each minimization step
in order to regain feasibility. This projection onto the constraints set is in
itself a non-trivial optimization problem and the need to solve it in every
iteration hinders the projected gradient methods and limits their efficiency
to only feasible sets that are \textquotedblleft simple to project
on\textquotedblright. Barrier or penalty methods likewise are based on
unconstrained optimization married with various \textquotedblleft
add-on\textquotedblright s that guarantee that the constraints are preserved.
Regularization methods embed the constraints into the objective function and
proceed with unconstrained solution methods for the \textquotedblleft
regularized\textquotedblright\ new objective function.

In contrast to these approaches, the superiorization methodology can be viewed
as an antipodal way of thinking. Instead of adapting unconstrained
minimization algorithms to handling constraints it adapts feasibility-seeking
algorithms to reduce target function values. This is done while retaining the
feasibility-seeking nature of the algorithm, and without paying a high
\textquotedblleft computational price\textquotedblright.

\textbf{Usefulness of the approach.} The usefulness of this approach relies on
two features: (i) \textbf{Computational}: feasibility-seeking is logically a
less-demanding task than seeking a constrained minimization point in a
feasible set. Therefore, letting efficient feasibility-seeking algorithms
\textquotedblleft lead\textquotedblright\ the algorithmic effort and modifying
them with inexpensive add-ons works well in practice. (ii)
\textbf{Applicational}: in some significant real-world applications the choice
of a target function is exogenous to the modeling and data collection which
give rise to the constraints. In such situations the limited confidence in the
usefulness of a chosen target function leads often to the recognition that,
from the application-at-hand point of view, there is no need, neither a
justification, to search for an exact constrained minimum. For obtaining
\textquotedblleft good results\textquotedblright,\ evaluated by how well they
serve the task of the application at hand, it is often enough to find a
feasible point that has reduced (not necessarily minimal) target function
value. In some operations research\ applications, the target functions are
costs or profits and are central to the model but in others the above
reasoning may still apply\footnote{Some support for this reasoning may be
borrowed from the American scientist and Noble-laureate Herbert Simon who was
in favor of \textquotedblleft satisficing\textquotedblright\ rather than
\textquotedblleft maximizing\textquotedblright. Satisficing is a
decision-making strategy that aims for a satisfactory or adequate result,
rather than the optimal solution. This is because aiming for the optimal
solution may necessitate needless expenditure of time, energy and resources.
The term \textquotedblleft satisfice\textquotedblright\ was coined by Herbert
Simon in 1956 \cite{simon}, see: https://en.wikipedia.org/wiki/Satisficing.}.

\textbf{Current research}. Current work on superiorization can be appreciated
from the materials on the Internet page \cite{bib-super-page}. In particular,
\cite{gth-sup4IA} and \cite{constantza15} are reviews of interest. Recent
research includes a variety of reports ranging from new applications in
industrial x-ray computed tomography \cite{schrapp14} to new mathematical
results on the foundation of superiorization such as strict Fej\'{e}r
monotonicity by superiorization of feasibility-seeking projection methods
\cite{cz14-feje}. A recent detailed description of previous work related to
superiorization can be found in \cite[Section 3]{cdhst14}.

\textbf{Linear superiorization.} Linear superiorization\textbf{ }(henceforth
abbreviated: LinSup) considers linear programming (LP) problems wherein the
constraints as well as the objective function are linear. The two questions
that we set out to explore experimentally here are: (i) Does LinSup provide a
feasible point whose target function value is lower than that obtained by
running the same feasibility-seeking algorithm without superiorization but
under otherwise identical conditions? and (ii) How does LinSup fare in
comparison with the Simplex method for solving LP problems? Based on our
computational experiments presented here, the answers to these two questions
are: \textquotedblleft yes\textquotedblright\ and \textquotedblleft very
well\textquotedblright, respectively.

An interesting and promising aspect of the current experiments is the
dependence of the results on the test problem sizes. We found that the
advantages of LinSup become monotonically more pronounced as the problem sizes
increase. We treated problems of up to $8,000$ linear inequalities and vectors
of up to $10,000$ components, but the trend is visible and if it persists
beyond these problem sizes then LinSup might well become a useful
computational tool for huge size problems. Admittedly, our preliminary work
presented here relies on randomly generated problems and these are not typical
of the kinds of problems that linear programming has been called to solve over
the years.

We show in Section \ref{sec:take-1} that LinSup finds a superior feasible
point, i.e., a feasible point with lower target function value. In Section
\ref{sec:take-2} we demonstrate the computational behavior of LinSup versus
the classical Simplex algorithm for linear optimization. The general framework
of superiorization appears in Section \ref{sec:frame} and LinSup is then
presented in Section \ref{sec:LS}. Our experimental results were generated
with MATLAB \cite{matlab} and are presented in Sections \ref{sec:take-1} and
\ref{sec:take-2}. We make concluding remarks in Section \ref{sec:conclusions}
and list a variety of questions for further research on LinSup. The Appendix
(Section \ref{sec:evolution}) briefly describes the technical changes and
modifications that the algorithmic structure of the superiorized version of a
basic algorithm has undergone in the published literature over the past
several years since it inception.

\section{The superiorization methodology\label{sec:frame}}

Consider a pair\ $(M,\mathcal{A})$ where $M,$ called a target set, is a given
subset of a given subset $Q$ of the $J$-dimensional Euclidean space $M\subset
Q\subseteq R^{J}$. Let $\mathcal{A}:Q\rightarrow Q$ be an algorithmic operator
that defines an iterative process, called the basic algorithm\footnote{It will
become truly an algorithm after a stopping rule will be added to it.},%
\begin{equation}
x^{0}\in Q,\text{ }x^{k+1}=\mathcal{A}(x^{k}),\text{ }k=1,2,\ldots
\label{eq:target-alg}%
\end{equation}
whose task is to find a point in the target set $M$. We, henceforth, refer to
such a pair $(M,\mathcal{A})$ as a superiorization pair. Let $\phi:Q\subseteq
R^{J}\rightarrow R$ be a given real-valued function, called a target function

The superiorization methodology is intended for constrained function reduction
problems of the following form.

\begin{problem}
\label{prob:sm}\textbf{The Constrained Function Reduction Problem}. Let
$(M,\mathcal{A})$ be a superiorization pair and let $\phi:Q\subseteq
R^{J}\rightarrow R$ be a target function. Find a point $x^{\ast}$ of $M$ whose
function $\phi$ value is less (but not necessarily minimal)\textbf{\textit{ }%
}than that of a point in $M$ that would have been reached by applying the
basic algorithm for finding a point of $M.$
\end{problem}

The superiorization methodology approaches this problem by investigating the
perturbation resilience of the basic algorithm, and then using proactively
such perturbations in order to \textquotedblleft force\textquotedblright\ the
perturbed algorithm obtained from the basic algorithm to do, in addition to
its original task, also target function reduction steps. The so perturbed
algorithm is called \textquotedblleft the superiorized version\ of the basic
algorithm\textquotedblright.

If the basic algorithm is computationally efficient and useful, in terms of an
application at hand, for finding a point of $M$ and if it is perturbation
resilient and the perturbations are simple and not expensive to calculate,
then the advantage of this method is that, for essentially the computational
cost of the basic algorithm, we are able to solve the constrained function
reduction problem by steering the iterates according to the target function
reduction perturbations. The superiorization methodology automatically
generates the superiorized version\ of the basic algorithm. The vector
$x^{\ast}$, obtained by applying the superiorized version\ of the basic
algorithm, need not be a minimizer of $\phi$ over $M.$ For further details
about the kinds of perturbation resilience that may be used consult, e.g.,
\cite[Definitions 4 and 9]{constantza15} or \cite{cdh10,cdhst14,hgdc12}.

The above definitions and terminology depend on what precise meaning we attach
to the statement \textquotedblleft Find a point $x^{\ast}$ of $M$%
\textquotedblright\ in Problem \ref{prob:sm}. In weak superiorization,
\textquotedblleft finding a point of $M$\textquotedblright\ is understood as
generating an infinite sequence $\{x^{k}\}_{k=0}^{\infty}$ that converges to a
point $x^{\ast}\in M,$ thus $M$ must be nonempty. In strong superiorization
\textquotedblleft finding a point of $M$\textquotedblright\ is understood as
finding a point $x^{\ast}$ that is $\varepsilon$-compatible with $M,$ for some
positive $\varepsilon,$ i.e., a point whose proximity function which measures
by how much it violates $M$ has value smaller or equal to $\varepsilon.$ Thus,
nonemptiness of $M$ need not be assumed. These notions were defined in
\cite{constantza15}.

Two significant special cases of superiorization pairs $(M,\mathcal{A})$ in
the above framework come to mind although other cases are also possible.

\begin{case}
\label{case:cfp}The target set $M$ is the solution set of a \textit{convex
feasibility problem} (CFP) of the form: Find a vector $x^{\ast}\in\cap
_{i=1}^{I}C_{i},$ where $C_{i}\subseteq R^{J}$ are closed convex subsets, thus
$M=\cap_{i=1}^{I}C_{i}.$ In this case the algorithmic operator and the basic
algorithm (\ref{eq:target-alg}) it entails can be any of the wide variety of
feasibility-seeking algorithms, see, e.g.,
\cite{bb96,bk13,CEG12,annotated,cccdh10,CZ97}.
\end{case}

\begin{case}
\label{case:2}The target set $M$ is the solution set of another constrained
minimization problem: $\mathrm{minimize}\left\{  f(x)\mid x\in\Omega\right\}
$ of an objective function $f$ over a feasible region $\Omega,$ thus
$M:=\{x^{\ast}\in\Omega\mid f(x^{\ast})\leq f(x)$ for all $x\in\Omega\}.$ In
this case the algorithmic operator and the basic algorithm
(\ref{eq:target-alg}) it entails can be any of the wide variety of constrained
minimization algorithms abundant in the literature.
\end{case}

In this paper we do linear superiorization which concentrates on the special
situation of Case \ref{case:cfp} wherein all constraint sets $C_{i}$ as well
as the target function $\phi$ are linear. Superiorization work with other
target functions such as total variation (TV) appears in, e.g.,
\cite{cdh10,cdhst14,hgdc12}. Superiorization work on Case \ref{case:2}, where
$M$ is the solution set of a maximum likelihood optimization problem appears
in \cite{gh13, wenma13, tie}.

\section{Linear superiorization\label{sec:LS}}

\subsection{The problem and the algorithm\label{subsec:linsup}}

Let the target set $M$ be%
\begin{equation}
M:=\{x\in R^{J}\mid Ax\leq b,\text{ }x\geq0\} \label{eq:linear-feas-prob}%
\end{equation}
where the $I\times J$ real matrix $A=(a_{j}^{i})_{i=1,j=1}^{I,J}$ and the
vector $b=(b_{i})_{i=1}^{I}\in R^{I}$ are given.

For a basic algorithm we pick a feasibility-seeking projection method.
Projections onto sets are used in many methods in optimization theory but here
projection methods refer to iterative algorithms that use projections onto
sets while relying on the general principle that when a family of, usually
closed and convex, sets is present, then projections onto the individual sets
are easier to perform than projections onto other sets (intersections, image
sets under some transformation, etc.) that are derived from the individual sets.

Projection methods may have different algorithmic structures, such as
block-iterative projections (BIP), see, e.g., \cite{dhc09,gordon} and
references therein, or string-averaging projections (SAP), see, e.g.,
\cite{cz-2-2014} and references therein, of which some are particularly
suitable for parallel computing, and they demonstrate nice convergence
properties and/or good initial behavior patterns. This class of algorithms has
witnessed great progress in recent years and its member algorithms have been
applied with success to many scientific, technological and mathematical
problems. See, e.g., the 1996 review \cite{bb96}, the recent annotated
bibliography of books and reviews \cite{annotated} and its references, the
excellent book \cite{CEG12}, or \cite{cccdh10}.

An important comment is in place here. A convex feasibility problem, mentioned
in Case \ref{case:cfp}, can be translated into an unconstrained minimization
of some proximity function that measures the feasibility violation of points.
For example, using a weighted sum of squares of the Euclidean distances to the
sets of the CFP as a proximity function and applying steepest descent to it
results in a simultaneous projections method for the CFP of the Cimmino type.
However, there is no proximity function that would yield the sequential
projections method for CFPs of the Kaczmarz type, see \cite{baillon12}.
Therefore, the study of feasibility-seeking algorithms for the CFP has
developed independently of minimization methods and it still vigorously does,
see the references mentioned above. Over the years researchers have tried to
harness projection methods for the convex feasibility problem to LP in more
than one way, see, e.g., Chinneck's book \cite{chinneck-book}. The mini-review
of relations between linear programming and feasibility-seeking algorithms in
\cite[Section 1]{nurminski15} sheds more light on this. Our results lead us to
wonder whether LinSup can serve such a cause.

The target function for linear superiorization will be%
\begin{equation}
\phi(x):=\left\langle c,x\right\rangle \label{eq:lin-objective}%
\end{equation}
where $\left\langle c,x\right\rangle $ is the inner product of $x$ and a given
$c\in R^{J}.$

In the footsteps of the general principles of the superiorization methodology,
as presented for general target functions $\phi$ in previous publications,
consult, e.g., the recent reviews \cite{gth-sup4IA} and \cite{constantza15},
we present the following linear superiorization algorithm. The input to the
algorithm consists of the problem data $A,$ $b,$ and $c$ of
(\ref{eq:linear-feas-prob}) and (\ref{eq:lin-objective}), respectively, a
user-chosen initialization point $\bar{y}$ and a kernel $0<\alpha<1$ (see item
1 in Subsection \ref{subsec:details}) with which the algorithm generates the
step sizes $\beta_{k,n},$ as well as an integer $N$ (see item 7 in Subsection
\ref{subsec:details}). All quantities in the algorithm that have not yet been
defined or explained are detailed in the Subsection \ref{subsec:details} below.

\begin{algorithm}
\label{alg_super}\textbf{Linear Superiorization (LinSup)}
\end{algorithm}

\begin{enumerate}
\item \textbf{set} $k\leftarrow0$

\item \textbf{set} $y^{k}\leftarrow\bar{y}$

\item \textbf{set }$\ell_{-1}\leftarrow0$

\item \textbf{while }stopping rule not met\textbf{ do}

\item $\qquad$\textbf{set} $n\leftarrow0$

\item \ \ \ \ \ \ \textbf{set} $\ell\leftarrow\operatorname{rand}(k,\ell
_{k-1})$

\item $\qquad$\textbf{set} $y^{k,n}\leftarrow y^{k}$

\item $\qquad$\textbf{while }$n$\textbf{$<$}$N$ \textbf{do}

\item $\qquad\qquad$\textbf{set} $\beta_{k,n}\leftarrow\eta_{\ell}$

\item \qquad\qquad\ \textbf{set} $z\leftarrow y^{k,n}-\beta_{k,n}%
\frac{\displaystyle c}{\displaystyle\left\Vert c\right\Vert _{2}}$

\item $\qquad\qquad$\textbf{set }$n\leftarrow n+1$

\item $\qquad\qquad$\textbf{set }$y^{k,n}\leftarrow z$

\item \ \ \ \ \ \ \ \ \ \ \ \ \textbf{set }$\ell\leftarrow\ell+1$

\item \ \ \ \ \ \ \textbf{end while}

\item \ \ \ \ \ \ \textbf{set }$\ell_{k}\leftarrow\ell$

\item \qquad\textbf{set }$y^{k+1}\leftarrow\mathcal{A}\left(  y^{k,N}\right)
$

\item \qquad\textbf{set }$k\leftarrow k+1$

\item \textbf{end while}
\end{enumerate}

\subsection{The Agmon-Motzkin-Schoenberg algorithm as the basic
algorithm\label{subsec:AMS}}

We use the projection method of Agmon-Motzkin-Schoenberg (AMS)
\cite{agmon,mot-schoen54}, see also, e.g., \cite[Algorithm 5.4.2]{CZ97}, as
the basic algorithm for feasibility-seeking represented by $\mathcal{A}$ in
step 16 of Algorithm \ref{alg_super}. Denote the half-spaces represented by
individual rows of (\ref{eq:linear-feas-prob}) by $H_{i},$%
\begin{equation}
H_{i}:=\{x\in R^{J}\mid\left\langle a^{i},x\right\rangle \leq b_{i}\},
\end{equation}
where $a^{i}\in R^{J}$ is the $i$-th row of $A$ and $b_{i}\in R$ is the $i$-th
component of $b$ in (\ref{eq:linear-feas-prob})$.$ The orthogonal projection
of an arbitrary point $z\in R^{J}$ onto $H_{i},$ has the closed-form%
\begin{equation}
P_{H_{i}}(z)=\left\{
\begin{array}
[c]{ll}%
z-\displaystyle\frac{\left\langle a^{i},z\right\rangle -b_{i}}{\Vert
a^{i}\Vert^{2}}a^{i}, & \text{if }\left\langle a^{i},z\right\rangle >b_{i},\\
z, & \text{if }\left\langle a^{i},z\right\rangle \leq b_{i}.
\end{array}
\right.  \label{eq:P(x)}%
\end{equation}

\begin{algorithm}
\textbf{The Relaxation Method of Agmon, Motzkin and Schoenberg (AMS)}{.
\label{algo:AMS}}

\textbf{Initialization}: $x^{0}\in R^{n}$ is arbitrary.

\textbf{Iterative step}: Given the current iteration vector $x^{k}$ the next
iterate is calculated by
\begin{equation}
x^{k+1}=\left\{
\begin{array}
[c]{ll}%
x^{k}-\lambda_{k}\displaystyle\frac{\left\langle a^{i(k)},x^{k}\right\rangle
-b_{i}}{\Vert a^{i(k)}\Vert^{2}}a^{i}, & \text{if }\left\langle a^{i(k)}%
,x^{k}\right\rangle >b_{i(k)},\\
x^{k}, & \text{if }\left\langle a^{i(k)},x^{k}\right\rangle \leq b_{i(k)}.
\end{array}
\right.
\end{equation}

\textbf{Relaxation parameters}: The parameters $\lambda_{k}$ are such that
$\epsilon_{1}\leq\lambda_{k}\leq2-\epsilon_{2},$ for all $k\;\geq\;0,$ with
some, arbitrarily small, $\epsilon_{1},\epsilon_{2}>0.$

\textbf{Control}: The control sequence $\{i(k)\}_{k=0}^{\infty}$ is almost
cyclic on $\{1,2,...,I\}$.
\end{algorithm}

This AMS cyclic feasibility-seeking algorithm goes cyclically through the
inequalities of (\ref{eq:linear-feas-prob}). To handle the nonnegativity
constraints in (\ref{eq:linear-feas-prob}) we just take the current iteration
vector in hand, after having done a full sweep of AMS through all $I$
row-inequalities, and set its negative components to zero while keeping the
others unchanged.

A corner stone of the superiorization methodology, in general as well as for
the linear case discussed here, is the perturbation resilience of the basic
algorithm that is used. The AMS algorithm is known to be bounded perturbation
resilience, this can be obtained from various previously published results,
see, e.g., \cite[Theorem 12]{cz12}, \cite{ndh12}.

\subsection{Implementation details and explanations\label{subsec:details}}

Here are the implementation details of our experimental work with the LinSup
Algorithm \ref{alg_super} presented in the next sections.

\begin{enumerate}
\item \textbf{Step-sizes of the perturbations. }The step sizes $\beta_{k,n}$
in Algorithm \ref{alg_super} must be such that $0<\beta_{k,n}\leq1$ in a way
that guarantees that they form a summable sequence $\sum_{k=0}^{\infty}%
\sum_{n=0}^{N-1}\beta_{k,n}<\infty,$ see, e.g., \cite{cz14-feje}. To this end
Algorithm \ref{alg_super} assumes that we have available a summable sequence
$\{\eta_{\ell}\}_{\ell=0}^{\infty}$ of positive real numbers generated by
$\eta_{\ell}=\alpha^{\ell}$ , where $0<\alpha<1$. Simultaneously with
generating the iterative sequence $\{y^{k}\}_{k=0}^{\infty},$ a subsequence
of\ $\{\eta_{\ell}\}_{\ell=0}^{\infty}$ is used to generate the step sizes
$\beta_{k,n}$ in step 9 of Algorithm \ref{alg_super}. The number $\alpha$ is
called the kernel of the sequence $\{\eta_{\ell}\}_{\ell=0}^{\infty}.$

\item \textbf{Controlling the decrease of the step-sizes of target function
reduction}. If during the application of Algorithm \ref{alg_super} the step
sizes $\beta_{k,n}$ decrease too fast then too little leverage is allocated to
the target function reduction activity that is interlaced into the
feasibility-seeking activity of the basic algorithm. This delicate balance can
be controlled by the choice of the index $\ell$ updates and separately by the
value of $\alpha$ whose powers $\alpha^{\ell}$ determine the step sizes
$\beta_{k,n}$ in step 9. In our work we adopt a strategy for updating the
index $\ell$ that was proposed and implemented for total variation (TV) image
reconstruction from projections by Prommegger and by Langthaler in \cite[page
38 and Table 7.1 on page 49]{prommegger} and in \cite{langthaler},
respectively. Instead of consecutively increasing $\ell$ by taking its value
as it was at the end of the last sweep of $N$ perturbations and starting the
new sweep from that last value, the Prommegger and Langthaler strategy
advocates to set $\ell$ at the beginning of every new iteration sweep (steps 5
and 6) to a random number between the current iteration index $k$ and the
value of $\ell$ from the last iteration sweep, i.e., $\ell_{k}%
=\operatorname{rand}(k,\ell_{k-1}).$ This strategy was denoted in those
reports by the name \textquotedblleft ATL2\textquotedblright\ and having
verified its utility for our work we adopted it in all our experiments. On the
other hand, the value of $\alpha$ with whose powers $\alpha^{\ell}$ the step
sizes $\beta_{k,n}$ are determined was experimented with and our computational
results in the sections below report on this for the experiments with LinSup
that we performed. Obviously, there are no such delicate balances in the
Simplex algorithm, although there were at the early stages (pivot strategies,
candidate list, etc.). Further work is needed to make LinSup more resistant to
the choice of parameters.

\item \textbf{No target function value comparisons}. Influenced by the results
of \cite{cz14-feje} we completely deleted from the algorithm the
decision-making test that compares the target function value at $z$ of step 10
in the perturbation inner-loop with the target function value at $y^{k}.$ This
decision-making test appeared in many previous formulations of the
superiorized version\ of the basic algorithm, see, e.g., step (xiv) of the
\textquotedblleft Superiorized Version of Algorithm $\mathbf{P}$%
\textquotedblright\ in \cite[page 5537]{hgdc12}. Since we were able to prove
our experimental claims without this test and since the mathematical treatment
in \cite{cz14-feje} also proceeded well without it we left it out. See also
the Appendix in Section \ref{sec:evolution} below.

\item \textbf{The proximity function}. To measure the feasibility-violation
(or level of agreement) of a point with respect to the target set $M$ we used
the following proximity function%
\begin{equation}
\Pr(x):=\frac{1}{2I}%
{\displaystyle\sum_{i=1}^{I}}
\frac{\left(  \left(  \left\langle a^{i},x\right\rangle -b_{i}\right)
_{+}\right)  ^{2}}{%
{\displaystyle\sum\limits_{j=1}^{J}}
\left(  a_{j}^{i}\right)  ^{2}}+\frac{1}{2J}%
{\displaystyle\sum\limits_{j=1}^{J}}
\left(  \left(  -x_{j}\right)  _{+}\right)  ^{2} \label{eq:prox}%
\end{equation}
where the plus notation means, for any real number $d,$ that $d_{+}%
:=\max(d,0).$ The proximity function is scaling invariant because it measures
exactly the weighted (with equal weights) sum of half square distances of the
point $x$ from all linear inequality constraints. These distances (see
(\ref{eq:P(x)})) are geometric entities insensitive to scaling.

\item \textbf{Initialization points of the algorithms}. In our experimental
studies the initialization point $\bar{y}$ was always chosen to be a
non-feasible point of the target set, $\bar{y}\notin M.$ Otherwise, there was
the danger that the algorithm would not move because of the AMS
feasibility-seeking step 16 of the algorithm. This was done in the following
manner. First $\bar{y}$ is randomly picked in the interval $[0,1]$ by the
algorithm or chosen otherwise by the user and its proximity function
(\ref{eq:prox}) to the set $M$ is calculated. If the proximity of $\bar{y}$ is
zero then we redefine a new $\bar{y}\leftarrow10\cdot\bar{y}$ and repeat these
$10$-fold increments until we found a point with nonzero proximity to serve as
the initialization point.

\item \textbf{Test problem generation}. We created target sets $M$ of various
sizes and a linear target function $\phi$ and run on them Algorithm
\ref{alg_super} with or without superiorization, as the case maybe, in our
experimental work reported below. Each problem of given $I\times J$ size was
created by defining a matrix $A$ whose elements are randomly chosen within the
interval $[-1,2]$ for all experiments. A vector $c$ was chosen randomly in
$[-2,3]$ always. To guarantee feasibility (nonemptiness) of the target set
$M,$ we defined $b$ by $b:=A\mathbf{1}+10\cdot\mathbf{1,}$ where $\mathbf{1}$
is the vector of all ones, and this guarantees that $\mathbf{1\in}M.$ In all
experiments we used problems with the following $I\times J$ sizes:
$80\times100,$ $200\times250,$ $400\times500,$ $800\times1,000,$
$2,000\times2,500,$ $4,000\times5,000,$ and $8,000\times10,000.$

We have not controlled the sparsity of the test problems and have not
investigated this issue. All we can say is that since all entries of the
matrices were uniformly distributed in the interval [-1,2], the probability of
any entry being equal to zero exactly is almost nonexistent. There are
millions of entries in total and many of them are probably close to zero but
overall it is accurate to say that the matrices that we generated were dense.

\item \textbf{The number }$N$\textbf{ of perturbation steps}. This number $N$
of perturbation steps that are performed prior to each application of the
feasibility-seeking operator $\mathcal{A}$ (in step 16) affects the
performance of the LinSup algorithm. It influences the balance between the
amounts of computations allocated to feasibility-seeking and those allocated
to target function reduction steps. A too large $N$ will make Algorithm
\ref{alg_super} spend too much resources on the perturbations that yield
target function reduction. In order to find an appropriate value of $N$ for
our work we created 10 problems of each of the problem sizes $80\times100,$
$200\times250,$ $400\times500,$ $800\times1,000,$ $2,000\times2,500,$ and 3
problems of the problem size $4,000\times5,000$. We applied Algorithm
\ref{alg_super} to each problem with the number $N$ being allowed to vary in
the range from $N=5$ to $N=100.$ All other parameters except $N$ were kept
equal in all runs (specifically, with kernel $\alpha=0.99,$ relaxation
parameters in the AMS algorithm in step 16 $\lambda_{k}=\lambda=1$ for all
$k\geq0,$ and initialization point $\bar{y}=10\cdot\mathbf{1}$ of appropriate
size for all problems.) The stopping rule for these experiments was when the
proximity function $\Pr(x)$ of (\ref{eq:prox}) dropped below the value
$\varepsilon=10^{-10}.$ We recorded for all runs the relative errors between
the linear target function value $\phi_{\text{\textrm{LinSup}}}$ obtained by
LinSup when it was stopped, and the linear objective function value
$\phi_{\text{Simplex}}$ obtained by the Simplex method when MATLAB reported
the solution has been reached, given by
\begin{equation}
RE:=\frac{\left\vert \phi_{\text{\textrm{LinSup}}}-\phi_{\text{Simplex}%
}\right\vert }{\left\vert \phi_{\text{Simplex}}\right\vert }%
.\label{eq:rel-error}%
\end{equation}
The table in Figure \ref{fig:N-table} contains averaged values of those
relative errors averaged over all problems of the same size. These data are
plotted in Figure \ref{fig:N-plot}. Based on these findings we decided to use
$N=30$ in all our subsequent computational experiments due to the observation
that, for all problem sizes, the decrease of relative error $RE$ became small
beyond this value of $N$. \qquad

\item \textbf{The relaxation parameters in the AMS algorithm}. In our
experiments we set all relaxation parameters in the AMS feasibility-seeking
algorithm represented by the algorithmic operator $\mathcal{A}$ and embodied
in step 16 of Algorithm \ref{alg_super} to $\lambda_{k}=\lambda=1.$ When doing
only feasibility-seeking proper with the AMS algorithm it has been frequently
shown in the literature that the relaxation parameters have significant effect
on the behavior of the algorithm, see, e.g., \cite[Subsections 11.2 and
11.5]{GTH}. However, here when the AMS algorithm is embedded in the LinSup
algorithm we observed that the relaxation parameters in the AMS algorithm have
a weak influence on the overall behavior and, therefore, they were set, at
this phase of the work, to $1$.

\item \textbf{Handling the nonnegativity constraints}. As mentioned above, the
nonnegativity constraints in (\ref{eq:linear-feas-prob}) are handled by taking
the current iteration vector in hand after having done a full sweep of AMS
through all $I$ row-inequalities of (\ref{eq:linear-feas-prob}) and setting
its negative components to zero while keeping the others unchanged.
\end{enumerate}

\section{Experimental Task 1: Linear superiorization finds a superior feasible
point\label{sec:take-1}}

Any of the large variety of projection methods to handle linear inequality
constraints feasibility-seeking can be used, but we choose for the basic
algorithm $\mathcal{A}$ the famous Agmon-Motzkin-Schoenberg (AMS) cyclic
feasibility-seeking projection method \cite{agmon,mot-schoen54}, known in the
image reconstruction literature as Algebraic Reconstruction Technique (ART)
for inequalities \cite[Subsection 11.2]{GTH}, see also \cite[Algorithm
5.4.2]{CZ97}.

Our aim in Task 1 is to experimentally validate or reject the following claim:

\begin{claim}
\label{claim:1}Consider two runs of the LinSup Algorithm \ref{alg_super} for
the same target set $M$ as in (\ref{eq:linear-feas-prob}), one with and the
other without superiorization. \textquotedblleft Without
superiorization\textquotedblright\ means that steps 5--15 in Algorithm
\ref{alg_super} are deleted and in step 16 one takes $y^{k,N}=y^{k}$ which
amounts to only running the feasibility-seeking basic algorithm $\mathcal{A}$
without any perturbations. Assume that other than that everything else is
equal in the two runs, such as the initialization point $\bar{y}$ and all
parameters associated with the application of the feasibility-seeking basic
algorithm $\mathcal{A}$ in step 16, as well as the stopping rule. Under these
circumstances the run \textquotedblleft with superiorization\textquotedblright%
\ \ will yield (i.e., stop at) a point $y^{\ast}$ whose $\phi(y^{\ast
}):=\left\langle c,y^{\ast}\right\rangle $ value will be smaller than
$\phi(y^{\ast\ast})$ of a point $y^{\ast\ast}$ at which the run
\textquotedblleft without superiorization\textquotedblright\ would stop.
\end{claim}

\begin{figure}
[ptb]
\begin{center}
\includegraphics[
natheight=5.333300in,
natwidth=14.750200in,
height=2.316in,
width=6.3581in
]%
{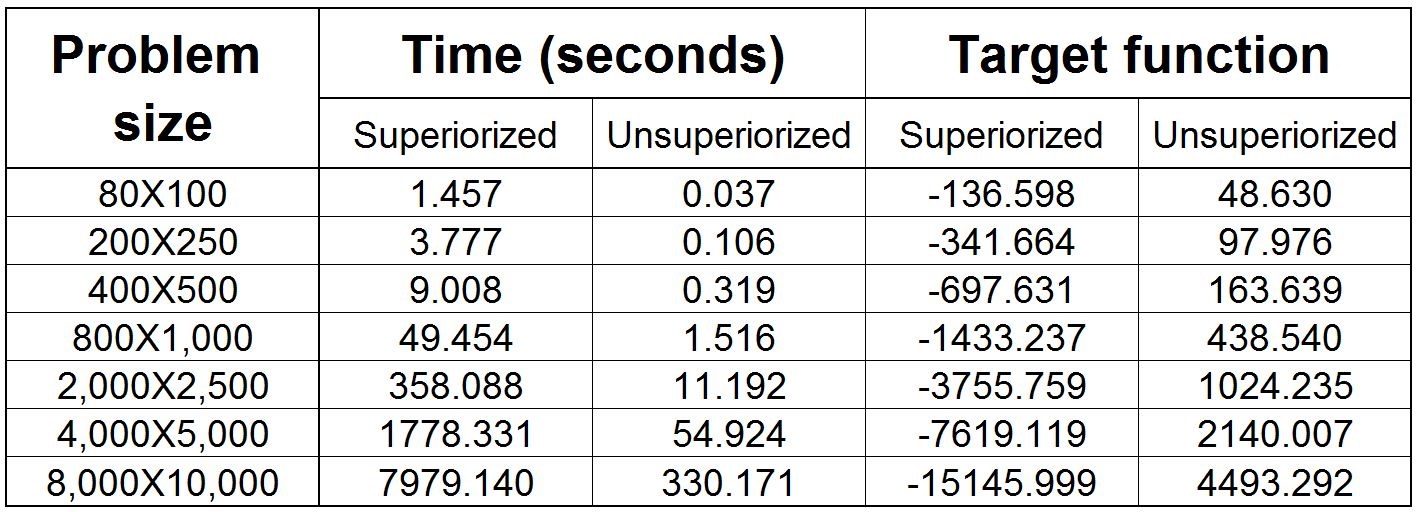}%
\caption{All values in this table are averaged over 10 different problems for
each problem size except for the last one ($8,000\times10,000$) for which we
made only one run. The execution times in seconds naturally show that
superiorization needs more time than plain feasibility-seeking. The two
right-hand side columns in the table confirm the truth of our Claim
\ref{claim:1} for the experiments that we performed.}%
\label{fig:task1-table}%
\end{center}
\end{figure}

\begin{figure}
[ptb]
\begin{center}
\includegraphics[
natheight=7.249700in,
natwidth=18.041700in,
height=2.3168in,
width=5.7251in
]%
{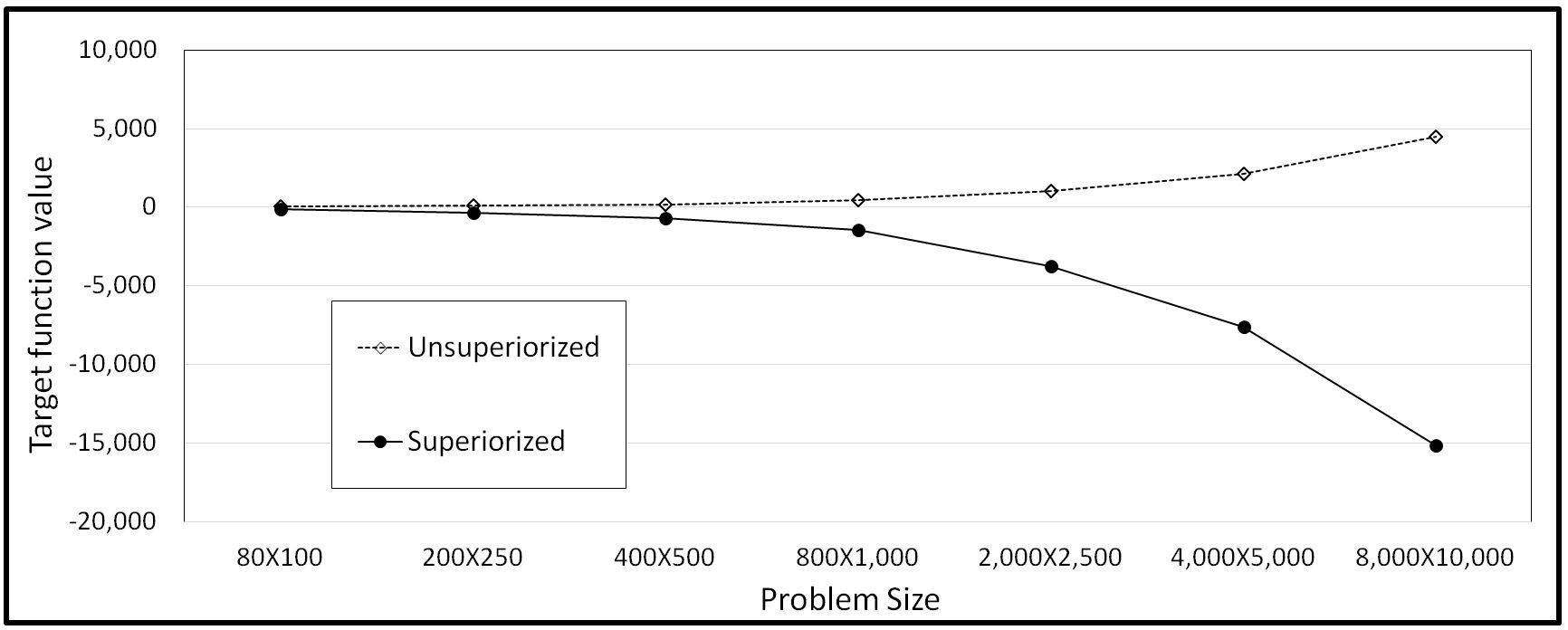}%
\caption{The data from the table in Figure \ref{fig:task1-table} is plotted
and shows that the gap between the target function values with and without
superiorization steadily increases with the increase in problem sizes.}%
\label{fig:task1-plot}%
\end{center}
\end{figure}

To prove this claim we created test problems as described in item 6 in
Subsection \ref{subsec:details}. On each such problem we ran LinSup without
superiorization and with superiorization and discovered that in all our
experiments Claim \ref{claim:1} is true. We ran all experiments with kernel
$\alpha=0.99,$ relaxation parameters in the AMS algorithm in step 16
$\lambda_{k}=\lambda=1$ for all $k\geq0,$ and initialization point $\bar
{y}=10\cdot\mathbf{1}$ of appropriate size for all problems. The stopping rule
for these experiments was when the proximity function $\Pr(x)$ of
(\ref{eq:prox}) dropped below the value $\varepsilon=10^{-20}.$ The number $N$
of perturbation steps that are performed prior to each application of the
feasibility-seeking operator $\mathcal{A}$ (in step 16) was, as decided in
item 7 in Subsection \ref{subsec:details}, $N=30.$ The execution times in
seconds, shown in the table in Figure \ref{fig:task1-table}, naturally show
that superiorization needs more time than plain feasibility-seeking. All
values in this table are averaged over 10 different problems for each problem
size except for the last one ($8,000\times10,000$) for which we made only one
run. The right-hand side columns in the table show the truth of our Claim
\ref{claim:1}. These data are plotted in Figure \ref{fig:task1-plot} and one
can clearly note that the trend persists and strengthens as the problem sizes increase.

Having generated our data as described in item 6 of Subsection
\ref{subsec:details}, the target function value actually depends also on the
size $J$ of the vector $x$. It is observed from the table in Figure
\ref{fig:task1-table} that when this size increases 10 times, the
corresponding target function values with and without superiorization both
roughly increase 10 times as well. From this point of view, the relative gap
between the target function values with and without superiorization is
consistent for different problem sizes.

\section{Experimental Task 2: Linear superiorization versus linear
optimization with the Simplex algorithm\label{sec:take-2}}

To compare the performance of LinSup with that of a linear optimization
algorithm we used MATLAB \cite{matlab} and chose the `Simplex' algorithm from
the `linprog' solver. We created test problems as described in item 6 in
Subsection \ref{subsec:details}. Since we wish to compare with the outputs and
execution times of the Simplex algorithm we first let MATLAB's Simplex
algorithm run on each test problem to ascertain that the `exitflag' that it
yields is `Function converged to a solution $x$'. If the test problem turned
out to be not solvable by the Simplex algorithm we discarded it in favor of
another test problem generated as described in item 6 in Subsection
\ref{subsec:details} for which Simplex outputs a solution.

Once a test problem was solved by Simplex we calculated the proximity $\Pr(x)$
of (\ref{eq:prox}) of the solution provided by the Simplex algorithm, which
was generally small. This proximity value was then used as the stopping rule
for the LinSup run on the same problem. When LinSup reached this proximity it
stopped and the iterate at stopping was its output solution.

Having forced the LinSup to run until it reached the same proximity as the
solution obtained by the Simplex algorithm, we recorded and compared the
target function values and the execution times for both. Based on the
experience gained in numerous experiments and runs we made decisions that
fixed all parameters except for one and report here on the performance of
LinSup and MATLAB's Simplex algorithm for several values of this parameter and
for different problem sizes. As said above, the number $N$ of perturbation
steps that are performed prior to each application of the feasibility-seeking
operator $\mathcal{A}$ was fixed to $N=30$ in all experiments. The
feasibility-seeking operator $\mathcal{A}$ (in step 16) was the AMS algorithm
of Subsection \ref{subsec:AMS} with fixed relaxation parameters $\lambda
_{k}=\lambda=1$ as in item 8 in Subsection \ref{subsec:details}. All other
implementation details were as in Subsection \ref{subsec:details}. We explored
the effect of different choices of the kernel $\alpha$ (in item 1 in
Subsection \ref{subsec:details}) on all runs.

All data presented in the following tables and plots is averaged over 10
different and independently-generated problems for each size from
$80\times100$ to $2,000\times2,500,$ 5 different and independently-generated
problems of size $4,000\times5,000$ and one problem of size $8,000\times
10,000.$

In addition to the relative error $RE$ of (\ref{eq:rel-error}) we recorded
here also the time ratio%
\begin{equation}
TR:=\frac{\text{execution time of LinSup}}{\text{execution time of Simplex}}.
\label{eq:time-ratio}%
\end{equation}

\subsection{The results and what they tell us\label{subsec:results-tell}}

\begin{figure}
[ptb]
\begin{center}
\includegraphics[
natheight=4.937200in,
natwidth=12.677300in,
height=2.3168in,
width=5.9049in
]%
{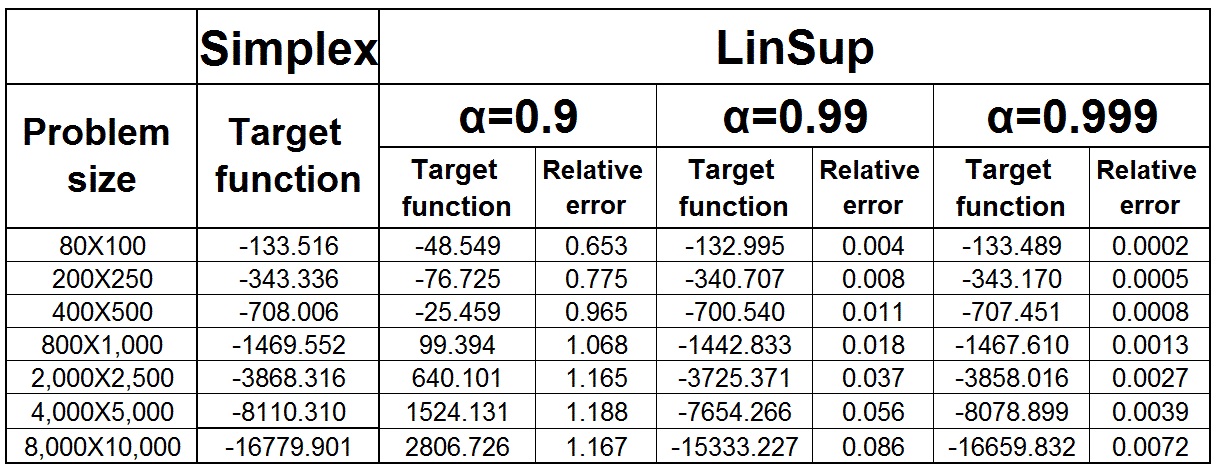}%
\caption{This table shows the target function $\phi(x)=\left\langle
c,x\right\rangle $ values for the Simplex algorithm alongside with the target
function values outputs by LinSup at stopping for 3 different values of the
kernel $\alpha.$ The relative errors $RE$ of (\ref{eq:rel-error}) are also
shown.}%
\label{fig:task2-target-table}%
\end{center}
\end{figure}

The table in Figure \ref{fig:task2-target-table} shows target function
$\phi(x)=\left\langle c,x\right\rangle $ values for the Simplex algorithm
alongside with the target function values outputs by LinSup at stopping for 3
different values of the kernel $\alpha.$ The relative errors $RE$ of
(\ref{eq:rel-error}) are also shown. With larger values of the kernel $\alpha$
its powers diminish slower, leaving more room for the target function
reduction perturbations to affect the outcome of LinSup. As the problem sizes
increase however the relative error $RE$ also increases.

\begin{figure}
[ptb]
\begin{center}
\includegraphics[
natheight=4.854200in,
natwidth=12.417000in,
height=2.3168in,
width=5.8842in
]%
{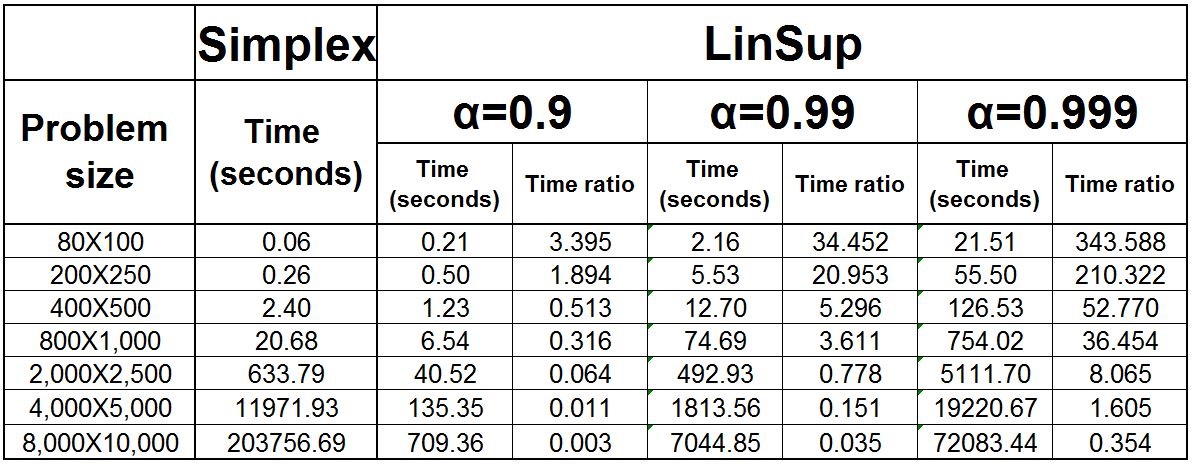}%
\caption{This table shows execution times in seconds for the Simplex algorithm
alongside with those of LinSup for 3 different values of the kernel $\alpha.$
The time ratios $TR$ of (\ref{eq:time-ratio}) are also shown.}%
\label{fig:linsup-task2-time-table}%
\end{center}
\end{figure}

The table in Figure \ref{fig:linsup-task2-time-table} shows execution times in
seconds for the Simplex algorithm alongside with those of LinSup for 3
different values of the kernel $\alpha.$ The time ratios $TR$ of
(\ref{eq:time-ratio}) are also shown. Here one observes that LinSup is fast
compared to the time of the Simplex algorithm.

\begin{figure}
[ptb]
\begin{center}
\includegraphics[
natheight=9.041600in,
natwidth=12.396200in,
height=2.316in,
width=3.1661in
]%
{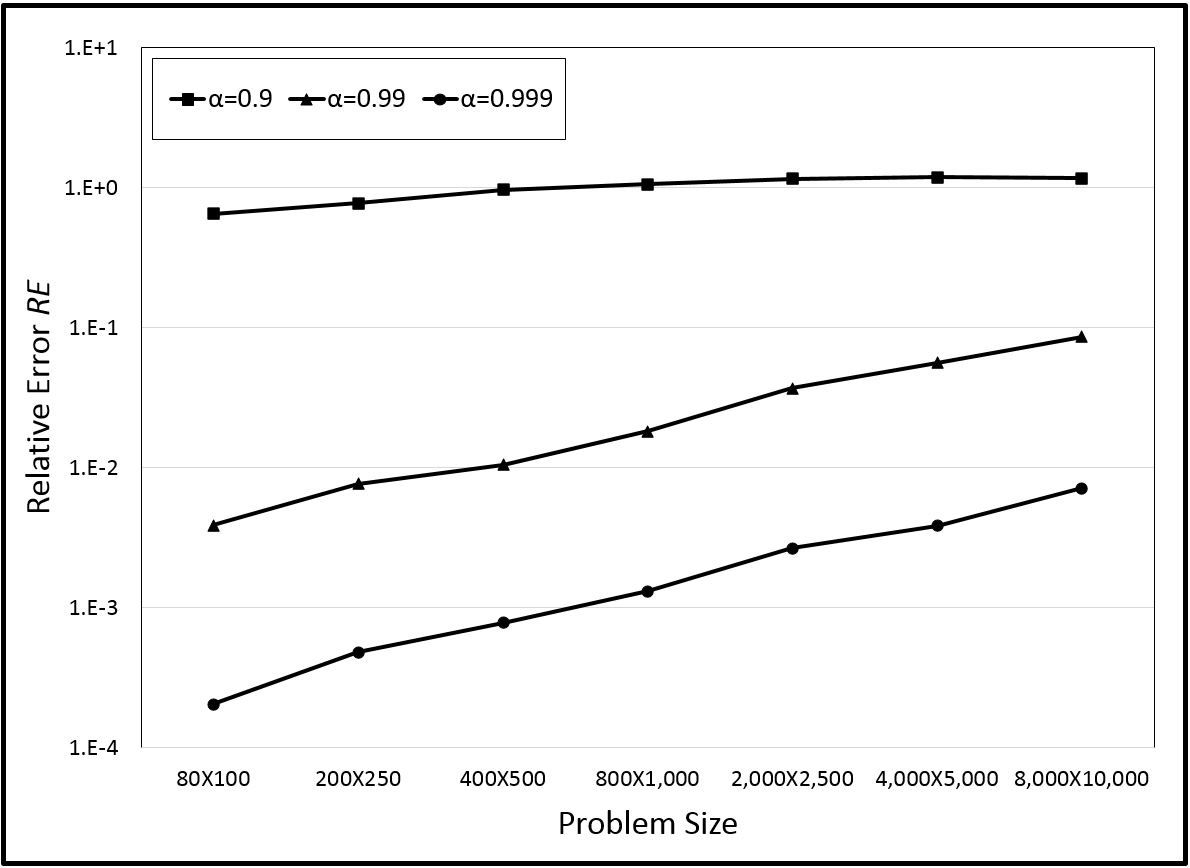}%
\caption{Plots of relative errors $RE,$ on a logarithmic scale, versus problem
sizes for LinSup with 3 different kernel $\alpha$ values based on the data
from the table in Figure \ref{fig:task2-target-table}. For each $\alpha$ the
relative error increases with the increase of problem sizes. For all problem
sizes the relative errors decrease with increasing value of $\alpha.$}%
\label{fig:gimel-task2-rel-err-log-plot}%
\end{center}
\end{figure}

\begin{figure}
[ptb]
\begin{center}
\includegraphics[
natheight=9.031200in,
natwidth=12.427400in,
height=2.3168in,
width=3.1782in
]%
{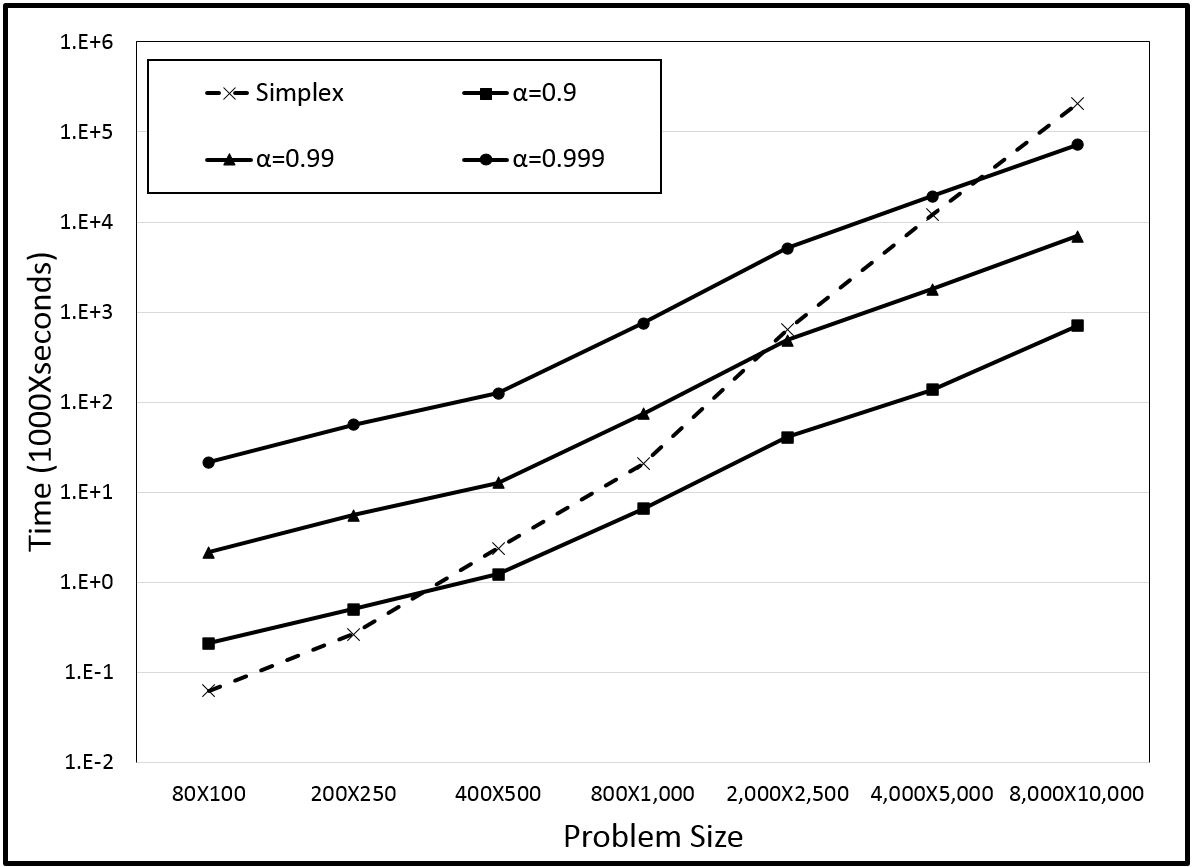}%
\caption{Execution times in thousands of seconds, on a logarithmic scale,
versus problem sizes of the Simplex algorithm and of LinSup for 3 kernel
$\alpha$ values. Observe the steep increase in time of the Simplex algorithm
(dashed line).}%
\label{fig:het-absolute-time-log-plot}%
\end{center}
\end{figure}

The Figures \ref{fig:gimel-task2-rel-err-log-plot}%
--\ref{fig:het-absolute-time-log-plot} are based on the data in the tables of
Figures \ref{fig:task2-target-table}--\ref{fig:linsup-task2-time-table}. Plots
of relative errors $RE$ versus problem sizes for LinSup with 3 different
kernel $\alpha$ values based on the data from the table in Figure
\ref{fig:task2-target-table} appear in Figure
\ref{fig:gimel-task2-rel-err-log-plot}. For each $\alpha$ the relative error
increases with the increase of problem sizes. For all problem sizes the
relative errors decrease with increasing value of $\alpha.$ For all problem
sizes the relative error is smaller for the larger value of $\alpha=0.999.$

\begin{figure}
[ptb]
\begin{center}
\includegraphics[
natheight=9.020900in,
natwidth=12.427400in,
height=2.3168in,
width=3.1816in
]%
{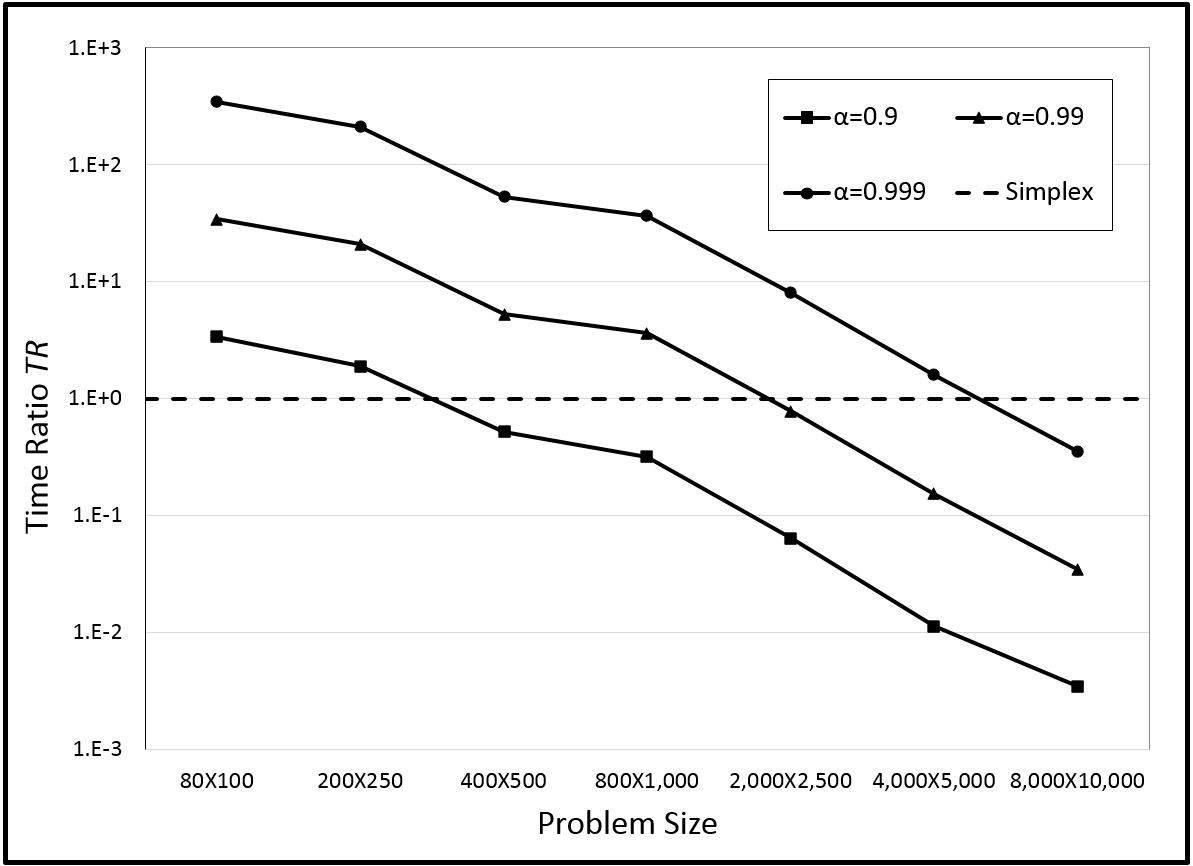}%
\caption{Plots of time ratios $TR,$ on a logarithmic scale, versus problem
sizes for LinSup with 3 different kernel $\alpha$ values based on the data
from the table in Figure \ref{fig:linsup-task2-time-table}. For each $\alpha$
the time ratio decreases with the increase of problem sizes. For all problem
sizes the time ratios decrease for decreasing value of $\alpha.$}%
\label{fig:dalett-task2-time-ratior-log-plot}%
\end{center}
\end{figure}

Plots of time ratios $TR$ versus problem sizes for LinSup with 3 different
kernel $\alpha$ values based on the data from the table in Figure
\ref{fig:linsup-task2-time-table} appear in Figure
\ref{fig:dalett-task2-time-ratior-log-plot}. For each $\alpha$ the time ratio
decreases with the increase of problem sizes. For all problem sizes the time
ratios decrease for decreasing value of $\alpha.$ This draws our attention to
the emerging conflict of choosing the kernel $\alpha.$ For better (smaller)
relative error choose it larger but for better (smaller) time ratio choose it
smaller. We see these trade-offs in the next figures as well.

\begin{figure}
[ptb]
\begin{center}
\includegraphics[
natheight=9.062400in,
natwidth=12.417000in,
height=2.3168in,
width=3.1644in
]%
{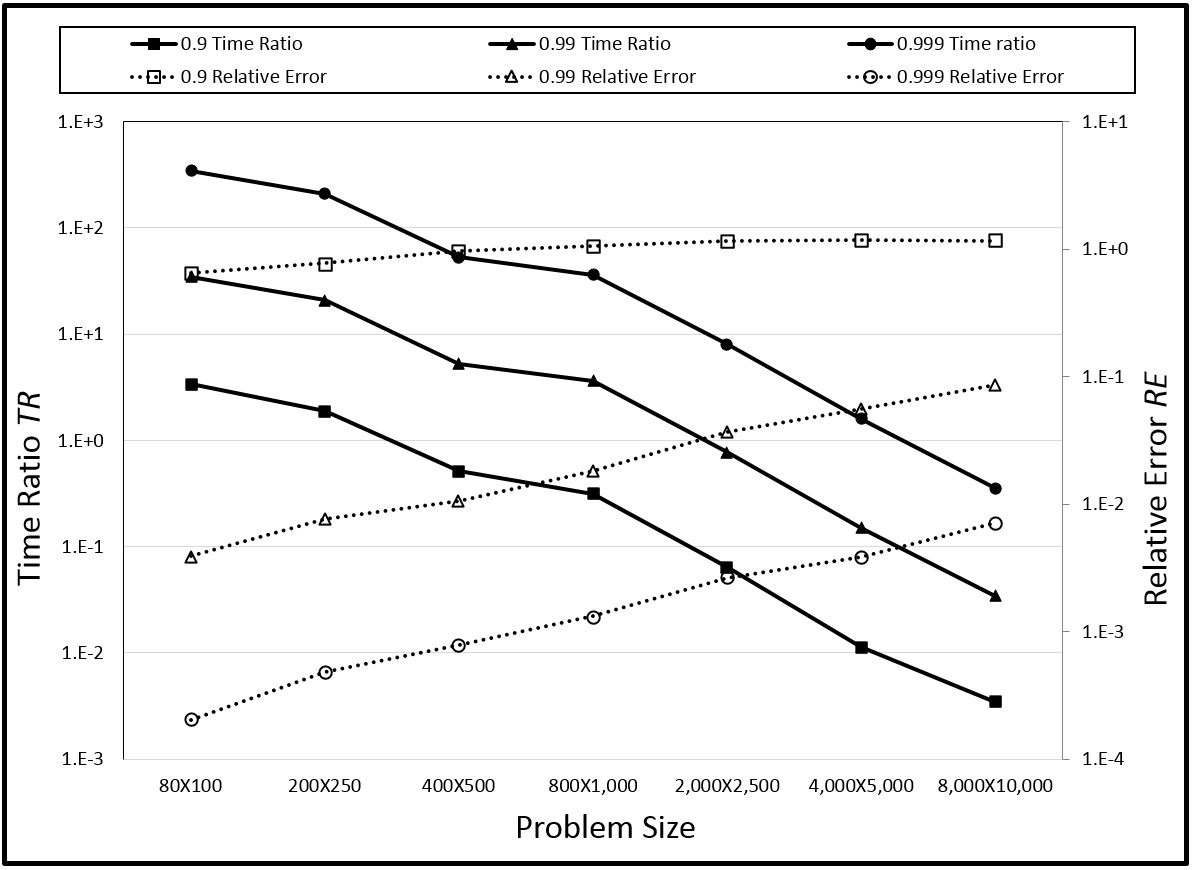}%
\caption{Here are the Figures \ref{fig:gimel-task2-rel-err-log-plot} and
\ref{fig:dalett-task2-time-ratior-log-plot} superimposed. This shows
graphically the trade-off between target function value reduction and speed in
the LinSup algorithm.}%
\label{fig:task2-target-and-time-superimposed-log}%
\end{center}
\end{figure}

Figure \ref{fig:hey-task2-target-absolute-plots} shows target function values
plotted against problem sizes for the 3 values of the kernel $\alpha.$ The
larger $\alpha=0.999$ allows for more resource investment of the LinSup
algorithm into function reduction steps. Thus, it yields target function
values that are closer to those obtained from the Simplex algorithm. Figure
\ref{fig:task2-target-and-time-superimposed-log} tells the story in a nutshell
by superimposing Figures \ref{fig:gimel-task2-rel-err-log-plot} and
\ref{fig:dalett-task2-time-ratior-log-plot}. This shows graphically the
trade-off between target function value reduction and speed in the LinSup algorithm.

Execution times in thousands of seconds versus problem sizes of the Simplex
algorithm and of LinSup for 3 kernel $\alpha$ values are depicted in Figure
\ref{fig:het-absolute-time-log-plot}. Observe the steep increase in time of
the Simplex algorithm (dashed line) for the larger sized problem. LinSup is
more moderate in the growth of needed execution times vis-a-vis the Simplex algorithm.

Our results show that there is a built-in \textquotedblleft
conflict\textquotedblright\ in choosing the parameters that govern the
delicate balance between the efforts that the LinSup algorithm invests in
feasibility-seeking and in function reduction with perturbations. But the
behavior of these results along increasing problem sizes leave room to hope
that with further increase of problem sizes LinSup will gain more ground and
become even a competitor to linear minimization algorithms. Observe that for
the problem in the last row of the tables in Figures
\ref{fig:task2-target-table} and \ref{fig:linsup-task2-time-table} LinSup with
$\alpha=0.999$ stops at target function value quite close to the one obtained
by the Simplex algorithm at about one third of the time it took the Simplex algorithm.

\subsection{Allowing the Simplex to terminate
suboptimally\label{subsec:suboptimal}}

LinSup is not intended to solve the LP problem but, as explained in Section
\ref{sec:frame}, to provide a feasible point with reduced (not necessarily
minimal) linear target function value. However, from the point of view of the
LP problem an output of LinSup can be considered a \textquotedblleft
reasonably good approximate solution of the LP problem\textquotedblright. This
raises the question, suggested by a referee, how would this compare with a
suboptimally terminated Simplex run. To take a preliminary look at this issue
we generated a $8,000\times10,000$ LP problem and let Simplex and LinSup run
on it. The LinSup was stopped when its iterates showed no further significant
changes (i.e., when $\frac{\displaystyle\left\Vert x^{k+1}-x^{k}\right\Vert
}{\displaystyle\left\Vert x^{k}\right\Vert }\leq10^{-16}$) and it was run with
two different values of $\alpha=0.99$ and $\alpha=0.995.$ The Simplex was not
allowed to run until optimality but stopped at a time that is just a little
longer than the time it took the LinSup runs to stop. These stopping decisions
enable us to compare LinSup with a suboptimally stopped Simplex on this
problem. Proximity function and linear target function calculation times after
each iteration were subtracted from the Simplex run times because they are not
an integral part of Simplex.

\begin{figure}
[ptb]
\begin{center}
\includegraphics[
natheight=9.031200in,
natwidth=19.999599in,
height=2.3168in,
width=5.0972in
]%
{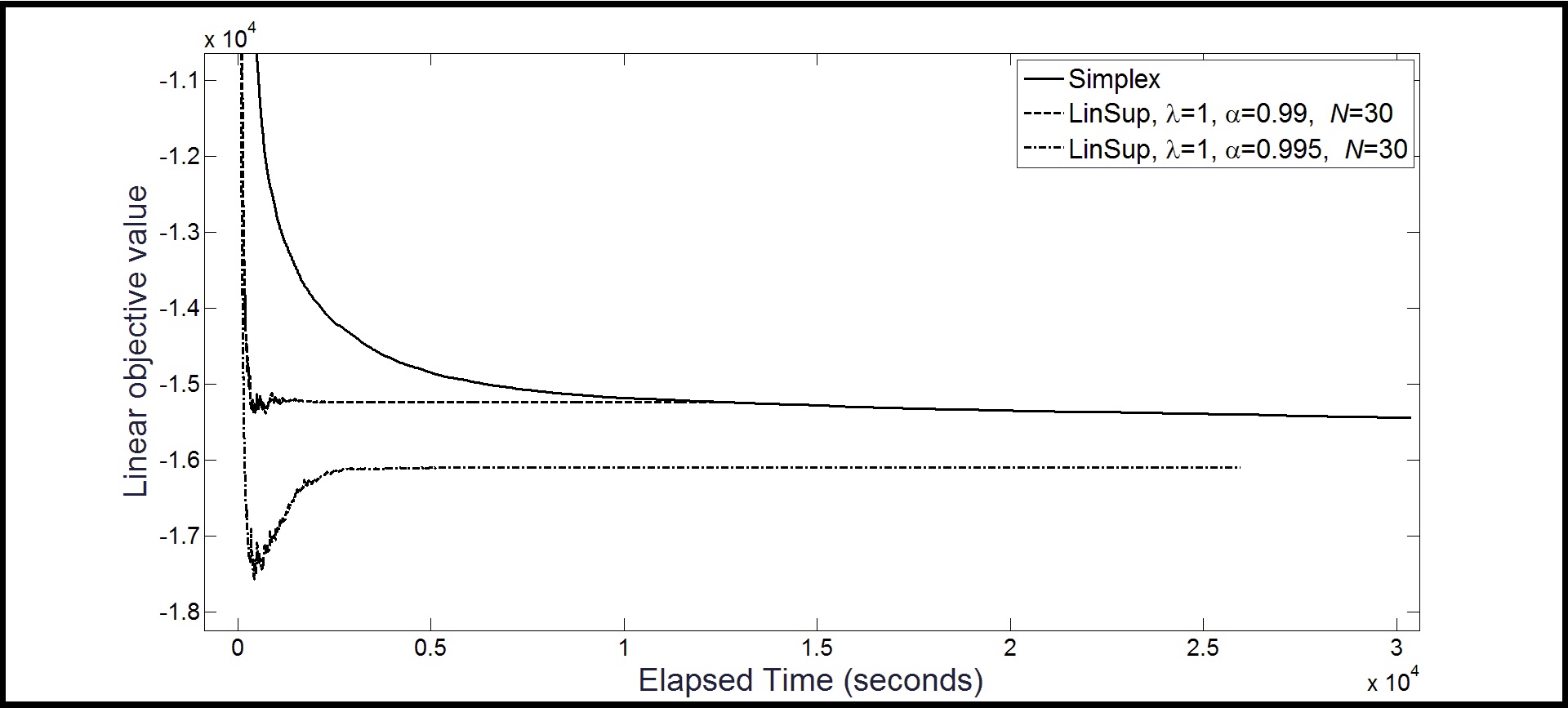}%
\caption{If this Simplex run would have been stopped suboptimally, say after
5,000 seconds, both runs of the LinSup would have yielded lower linear target
function values.}%
\label{fig:new-objective}%
\end{center}
\end{figure}

\begin{figure}
[ptb]
\begin{center}
\includegraphics[
natheight=9.031200in,
natwidth=19.999599in,
height=2.3168in,
width=5.0972in
]%
{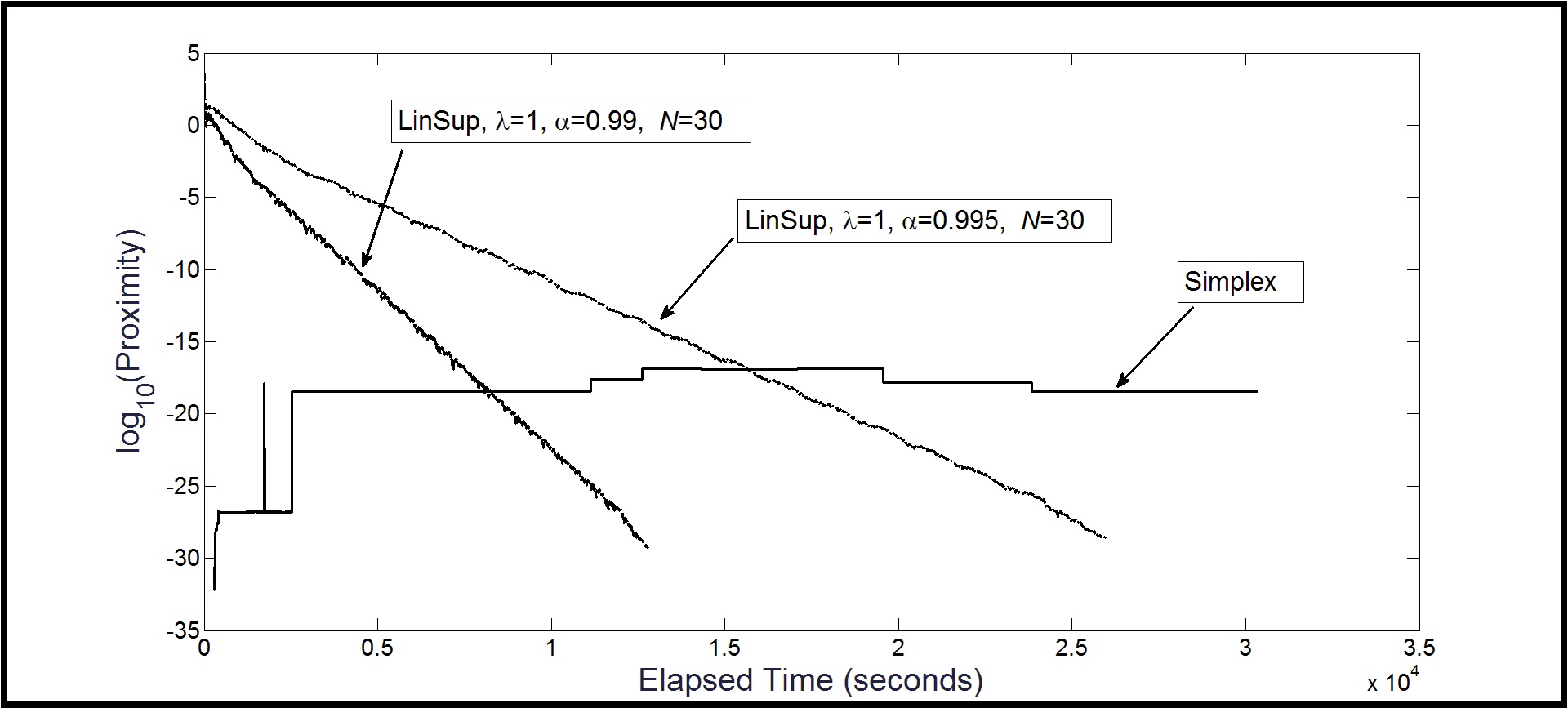}%
\caption{After 20,000 seconds, both LinSup runs would have a lower proximity
than the Simplex and the one with the higher kernel value $\alpha$ would even
have a lower linear target function value.}%
\label{fig:new-prox}%
\end{center}
\end{figure}

These results are depicted in Figures \ref{fig:new-objective} and
\ref{fig:new-prox}. Although far from being fully explored, the results show
that if this Simplex run would have been stopped suboptimally, say after 5,000
seconds, both runs of the LinSup would have yielded lower linear target
function values, as seen in Figure \ref{fig:new-objective}. At this point in
time Simplex would have delivered an output with better feasibility, i.e.,
lower proximity value. However, at a later point in time, say after 20,000
seconds, both LinSup runs would have a lower proximity than the Simplex as
seen in Figure \ref{fig:new-prox} and the one with the higher kernel value
$\alpha$ would even have a lower linear target function value.

This hints at the possible advantages of LinSup for large LP problems. Looking
at the output of LinSup as a \textquotedblleft reasonably good approximate
solution of the LP problem\textquotedblright, LinSup not only converges to
such a solution faster than it takes Simplex to solve a problem to machine
precision accuracy, but also faster than a suboptimally stopped Simplex.
Admittedly, this and the other experiments presented here call for further
work, see Section \ref{sec:conclusions}.

\section{Conclusions\label{sec:conclusions}}

Linear superiorization (LinSup) is not, as far as we know at this time, a
minimization method. Finding a constrained minimum point with it cannot be
guaranteed. What it does is to steer feasibility-seeking algorithms toward
points with lesser (not necessarily minimal) linear target function values.
The computationally-efficient feasibility-seeking algorithms that use
projections onto the convex closed sets of the constraints, embodied in
LinSup, are particularly successful for the linear case. The perturbations to
reduce the linear target function values need no effort other than using $-c$
as a direction of descent. Therefore, previous work on the superiorization
methodology in general (see the references mentioned in the Introduction and
in the Appendix) along with the proof of concept experimental work presented
here suggest that LinSup is potentially a viable option to handle large LP problems.

Our results show that LinSup indeed finds a superior feasible point. That the
increase in execution times as function of problem sizes of LinSup is more
moderate than that of the Simplex algorithm. This motivates us in formulating
the following conjecture.

\begin{conjecture}
There exists a size-level of huge problems above which LinSup will perform
better than linear minimization algorithms. Maybe that this will call for
using feasibility-seeking projection methods inside LinSup that lend
themselves to parallelization, such as block-iterative projections (BIP) or
string-averaging projections (SAP) methods mentioned in Subsection
\ref{subsec:linsup} above.
\end{conjecture}

Many questions present themselves for further research based on the current
work. Here is a telegraphic list of some potentially interesting directions:

(i) Expand the computational work to larger problem sizes and differently
generated problems.

(ii) Test LinSup on a larger class of test problems than those used here such
as LP benchmark test problems from Netlib (http://www.netlib.org/) or other repositories.

(iii) Study LinSup with additional feasibility-seeking projection methods that
lend themselves to parallelization, such as block-iterative projections (BIP)
or string-averaging projections (SAP) methods.

(iv) Investigate the parameters' effects on the behavior of LinSup by
repeating experiments with different values of: The number $N$ of perturbation
steps that are performed prior to each application of the feasibility-seeking
operator $\mathcal{A}$, the relaxation parameters $\lambda_{k}$ in the
feasibility-seeking embedded basic algorithm, the kernel $\alpha$ with which
the step-sizes $\beta_{k,n}$ are generated.

(v) Advance the mathematical analysis of LinSup.

(vi) Repeat the above comparisons for additional linear optimization
algorithms such as `interior-point' or `active-set' in MATLAB or others.

(vii) Investigate the inconsistent case wherein the target set $M$ of
(\ref{eq:linear-feas-prob}) is empty and is replaced, e.g., by the set of
closest points to all constraints according to some proximity function. Linear
programming algorithms might not work but LinSup can still furnish a useful result.

(viii) Study LinSup for sparse linear constraints for which some projection
methods have already demonstrated their effectiveness as feasibility-seeking algorithms.

\section{Appendix: The algorithmic evolution of
superiorization\label{sec:evolution}}

The algorithmic structure of the superiorized version of a basic algorithm has
undergone changes and modifications over the past several years since it
inception. All changes preserve the underlying basic methodology and it is
useful to briefly review them here. In \cite{bdhk07} superiorization appeared
although the words superiorization and perturbation resilience were not yet in
use there. It built on some earlier theoretical work in \cite{brz06,brz08}.
The pseudocode on the right-hand side column of page 543 in \cite{bdhk07}
constitutes the first superiorization algorithm. The step sizes $\beta$ there
(line 9) are simply halved and there is one function reduction step (line 3)
for each sweep of the feasibility-seeking algorithm $\mathbf{P}$ (in line 6).
There are two decision-making steps (on lines 4 and 7). This algorithm was
used in Scott Penfold's thesis work \cite{penfold-thesis} (see also
\cite{scott-book}) and in the paper \cite{pscr10}. The Algorithms 2 and 3
(TVS1-DROP and TVS2-DROP, respectively) in \cite{pscr10} relate to different
variants of the built-in feasibility-seeking algorithm DROP of \cite{DROP},
not to different superiorization methods. In the same paper \cite{pscr10} the
expensive decision-making step (line 12) in both algorithms was removed
without adverse effects thus allowing significant time savings.

In the \textquotedblleft Superiorized version of algorithm $\mathbf{P}%
$\textquotedblright\ on page 6 of \cite{cdh10} both decision-making steps from
earlier versions still appear, this time in a single line (line xiii).
However, in the subsequent \cite{hgdc12} the expensive decision-making step of
line 12 in \cite{pscr10} is not to be seen anymore and, additionally, the
negative subgradient (negative gradient if the function is differentiable) is
replaced by any direction of \textquotedblleft non-ascend\textquotedblright%
\ for the function $\phi$ that is superiorized. This last variation is not
needed for total variation (TV) superiorization but might come to use for
other functions $\phi$. Another new ingredient in \cite{hgdc12} which is very
useful is the ability to do in an inner loop (from line vii till line xvii)
$N$ steps of function reduction ($N$ is user-determined) for each sweep of the
feasibility-seeking algorithm, denoted by $\mathbf{P}_{T}$ there (in line xviii).

An interesting comparative study appears in \cite{cdhst14}. The algorithm
there is called \textquotedblleft Superiorized Version of the Basic
Algorithm\textquotedblright\ and appears on pp. 737--738 wherein the
feasibility-seeking algorithm is denoted by $A_{C}$ (in step 18).

In \cite[Algorithm 4.1]{cz14-feje} another step forward was made by (i)
allowing the number $N$ (of \cite{hgdc12}) to vary from one iteration to
another so that $N$ is replaced by $N_{k}$ where $k$ is the iteration index,
and (ii) discarding the second decision-making check that was on line 14 in
\cite{cdhst14} and in earlier algorithms for superiorization. The replacement
of $N$ by $N_{k}$ is mathematically valid but, to the best of our knowledge,
has not yet been experimented with by anyone.

Two important recent works on implementations of the superiorization algorithm
appear in \cite{langthaler} and \cite{prommegger}. One important additional
modification in those, that we adopted in our work (see item 2 in Subsection
\ref{subsec:details} above), is the way of controlling the perturbations'
step-sizes $\beta_{k,n}$ in Algorithm \ref{alg_super} via a special strategy
of updating the index $\ell$ at each sweep of iterations.\bigskip

\textbf{Acknowledgments.} We gratefully acknowledge some preliminary
discussions with Ran Davidi and John Chinneck. We thank the two anonymous
referees for their constructive comments which helped us improve the paper.
The MATLAB programming was skillfully performed with great enthusiasm and
devotion by Yehuda Zur, for which we are indebted to him. This work was
supported by Research Grant No. 2013003 of the United States-Israel Binational
Science Foundation (BSF) and by Award No. 1P20183640-01A1 of the National
Cancer Institute (NCI) of the National Institutes of Health (NIH).\bigskip

\textbf{Comment}. Final version preprints of the author's papers cited in the
references list below are available at:
http://math.haifa.ac.il/yair/censor-recent-pubs.html. Other papers on
superiorization cited below have their abstracts and DOI codes posted on:
http://math.haifa.ac.il/yair/bib-\newline superiorization-censor.html.

%

\begin{figure}
[ptb]
\begin{center}
\includegraphics[
natheight=6.833800in,
natwidth=13.989200in,
height=2.0833in,
width=4.2367in
]%
{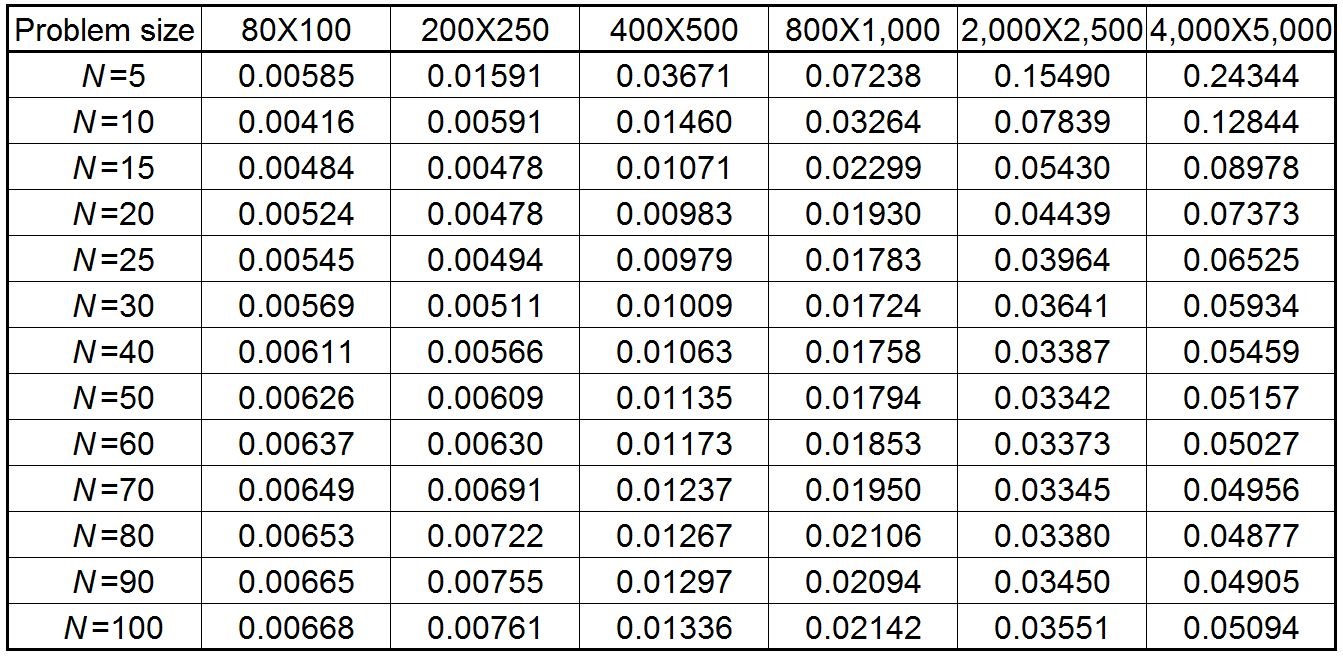}%
\caption{The realtive errors $\left\vert \phi_{\text{\textrm{LinSup}}}%
-\phi_{\text{Simplex}}\right\vert /\left\vert \phi_{\text{Simplex}}\right\vert
$ between the linear target function value $\phi_{\text{\textrm{LinSup}}}$
obtained by LinSup when it was stopped and the linear objective function value
$\phi_{\text{Simplex}}$ obtained by the Simplex method when MATLAB reported
the solution has been reached, for different values of the number $N$ of step
8 of Algorithm \ref{alg_super}. The numbers in the table are average values
over several problems of each size (see the text).}%
\label{fig:N-table}%
\end{center}
\end{figure}
%

\begin{figure}
[ptb]
\begin{center}
\includegraphics[
natheight=7.542000in,
natwidth=17.521099in,
height=2.316in,
width=5.3454in
]%
{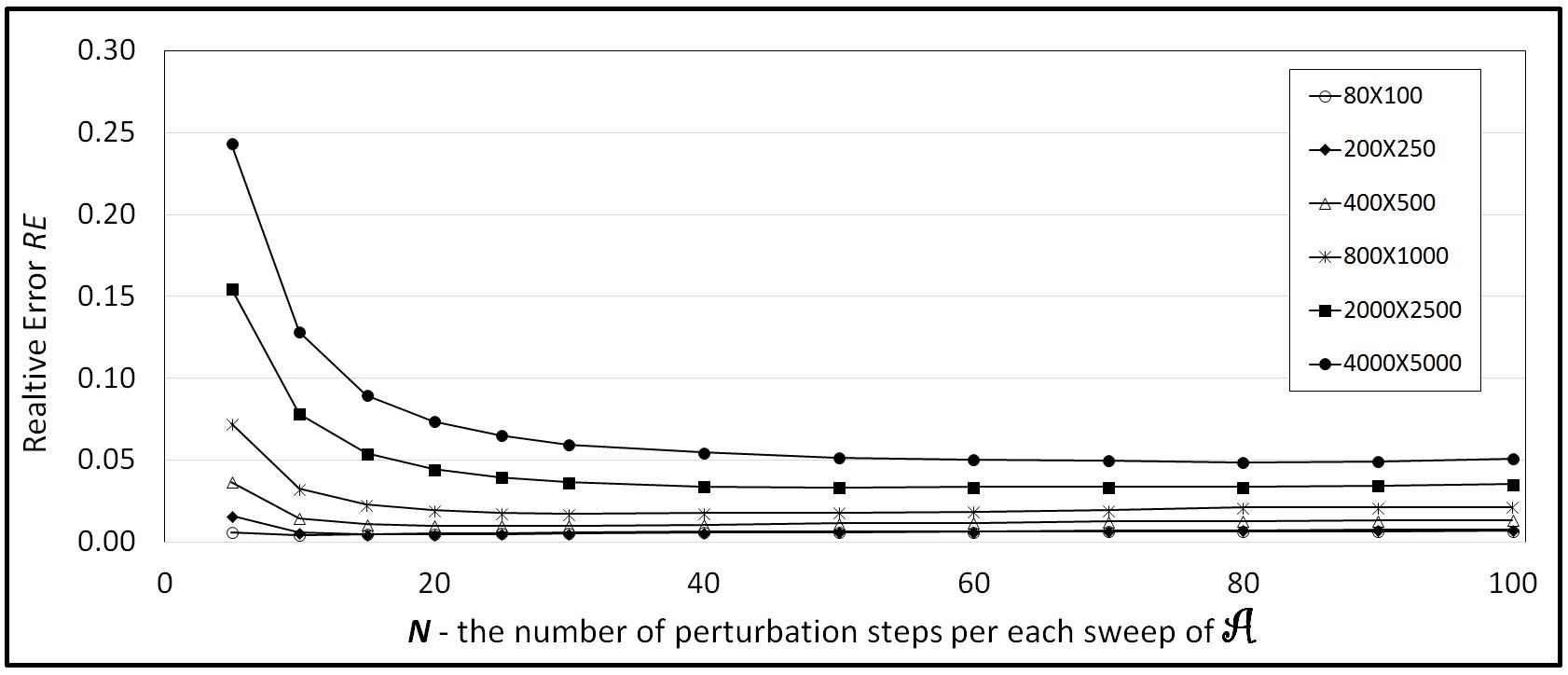}%
\caption{Plots of the data included in the table of Figure \ref{fig:N-table}.
Based on these findings we decided to use $N=30$ in all our subsequent
computational experiments.}%
\label{fig:N-plot}%
\end{center}
\end{figure}
%

%

%

%

%

%

%

\begin{figure}
[ptb]
\begin{center}
\includegraphics[
natheight=8.708700in,
natwidth=13.301700in,
height=2.316in,
width=3.5232in
]%
{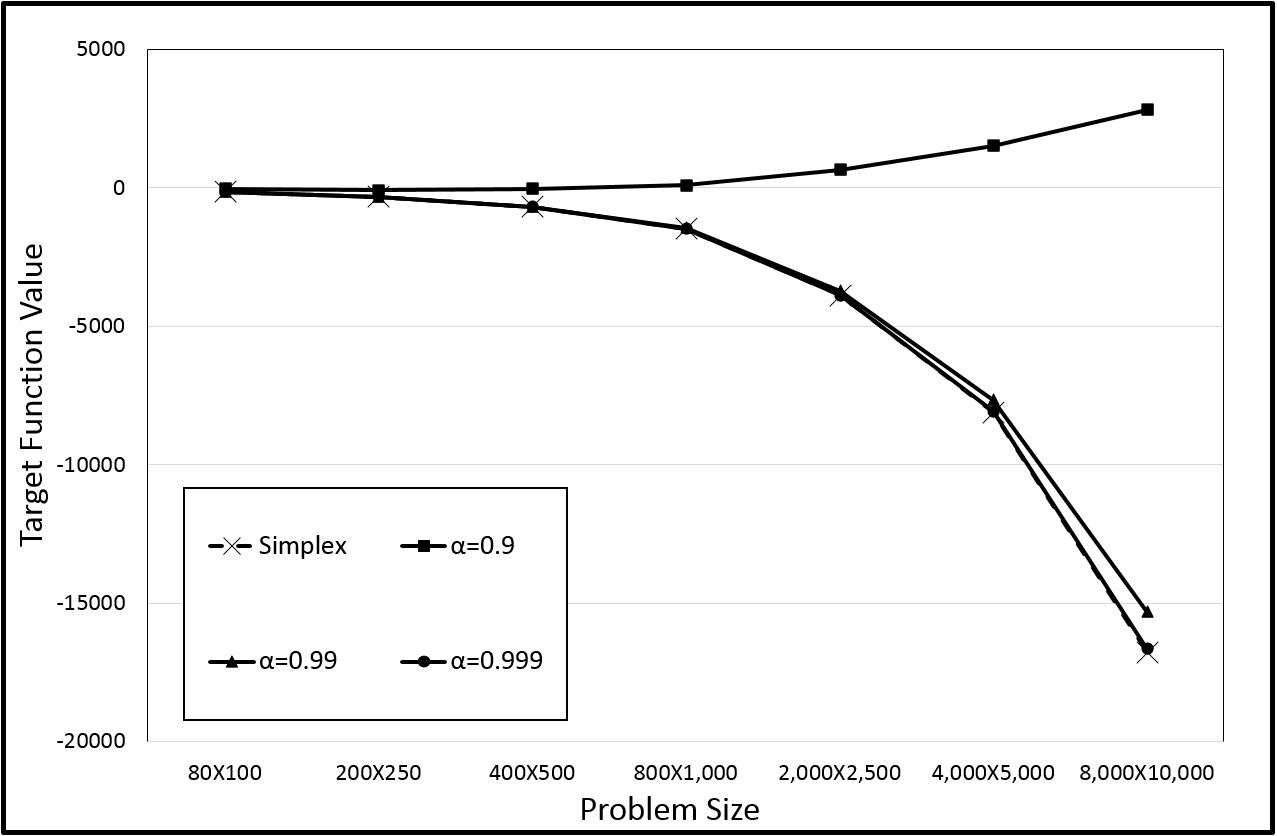}%
\caption{Target function values plotted against problem sizes for the 3 values
of the kernel $\alpha.$ The larger $\alpha=0.999$ allows more resource
investment of the LinSup algorithm into function reduction steps. thus, yields
target function values that are close to those obtained from the Simplex
algorithm.}%
\label{fig:hey-task2-target-absolute-plots}%
\end{center}
\end{figure}
%

%

%

%

\end{document}